\documentclass[letterpaper, 10 pt, conference]{ieeeconf}% Comment this line out

\usepackage{times}  %Required
\usepackage{helvet}  %Required
\usepackage{courier}  %Required
\usepackage{url}  %Required
\usepackage{graphicx}  %Required

\usepackage{algorithm}
\usepackage{algorithmic}

\IEEEoverridecommandlockouts                              % This command is only
\usepackage{fullpage,enumerate,setspace,changepage,multirow,rotating,xcolor}
\usepackage{graphicx,amsmath,amssymb,amsfonts,url, algorithm, algorithmic, etoolbox}
\usepackage{array,subfigure}
\usepackage[small]{caption}
\AtBeginEnvironment{align}{\setcounter{subeqn}{0}}% Reset subequation number at start of align
\newcounter{subeqn} \renewcommand{\thesubeqn}{\theequation\alph{subeqn}}%
\newcommand{\subeqn}{%
	\refstepcounter{subeqn}% Step subequation number
	\tag{\thesubeqn}% Label equation
}
\def\cP{\mathcal{P}}

\def\smskip{\smallskip}

\def\texitem#1{\par\smskip\noindent\hangindent 25pt
               \hbox to 25pt {\hss #1 ~}\ignorespaces}

\newcommand{\BEAS}{\begin{eqnarray*}}
\newcommand{\EEAS}{\end{eqnarray*}}
\newcommand{\BEA}{\begin{eqnarray}}
\newcommand{\EEA}{\end{eqnarray}}
\newcommand{\BEQ}{\begin{eqnarray}}
\newcommand{\EEQ}{\end{eqnarray}}
\newcommand{\BIT}{\begin{itemize}}
\newcommand{\EIT}{\end{itemize}}
\newcommand{\BNUM}{\begin{enumerate}}
\newcommand{\ENUM}{\end{enumerate}}

\newcommand{\BA}{\begin{array}}
\newcommand{\EA}{\end{array}}

\newcommand{\reals}{\mathbb{R}}

\newif\ifpagenumbering
\pagenumberingtrue

\pagenumberingfalse

\newsavebox{\theorembox}
\newsavebox{\lemmabox}
\newsavebox{\remarkbox}
\newsavebox{\assbox}
\savebox{\theorembox}{\noindent\bf Theorem}
\savebox{\lemmabox}{\noindent\bf Lemma}
\savebox{\remarkbox}{\noindent\bf Remark}
\savebox{\assbox}{\noindent\bf Assumption}

\newtheorem{lemma}{Lemma}

\makeatletter
\newcommand{\thickhline}{%
	\noalign {\ifnum 0=`}\fi \hrule height 1pt
	\futurelet \reserved@a \@xhline
}

\title{\LARGE \bf Sequential Chance Optimization For Flow-Tube Based Control Of Probabilistic Nonlinear Systems}
\author{Ashkan M. Jasour, Brian C. Williams\\
	MIT, Computer Science and Artificial Intelligence Laboratory \\
	\{jasour,williams@csail.mit.edu\} 
\thanks{
	This work was supported by Boeing grant MIT-BA-GTA-1.}
}

\begin{document}

\onecolumn

\maketitle
\thispagestyle{empty}
\pagestyle{empty}

%%%%%%%%%%%%%%%%%%%%%%%%%%%%%%%%%%%%%%%%%%%%%%%%%%%%%%%%%%%%%%%%%%%%%%%%%%%%%%%%%%%%%%%%%%%%%%%%%%%%%%%%%%%%%%%%%%%%%%%%%%%%%%%%

\begin{abstract}
	
In this paper, we address the problem of closed-loop control of nonlinear dynamical systems subjected to probabilistic uncertainties. More precisely, we design time-varying polynomial feedback controllers to follow the given nominal trajectory and also, for safety purposes, remain in the tube around the nominal trajectory, despite all uncertainties. We formulate this problem as a chance optimization problem where we maximize the probability of achieving control objectives. To address control problems with long planning horizons, we formulate the single large chance optimization problem as a sequence of smaller chance optimization problems. To solve the obtained chance optimization problems, we leverage the theory of measures and moments and obtain convex relaxations in the form of semidefinite programs. We provide numerical examples on stabilizing controller design and motion planning of uncertain nonlinear systems to illustrate the performance of the proposed approach.

% Provided approach deals with both bounded and unbounded probabilistic uncertainties and also long planning horizons. 
	 
\end{abstract}
%\begin{IEEEkeywords}
%	measure and moment theory, polynomial systems, sum of squares polynomials, semialgebraic set, SDP relaxation.
%\end{IEEEkeywords}
%**********************************************
\section{Introduction}

In this paper, we address the problem of nonlinear closed-loop controller design for nonlinear dynamical systems subjected to probabilistic uncertainties. This problem has many applications in control and robotics such as motion planning of autonomous systems in the presence of noise and disturbances and stability problem of uncertain dynamical systems. The problem of controller design in the presence of probabilistic uncertainties is hard and challenging. Because it requires computation and control of time evolution of
the probability distribution of states of the system to meet the control objectives.

There are different techniques to deal with uncertainties in the nonlinear controller design problem. Robust approaches consider bounded uncertainty sets and design controllers that are valid for all possible realization of the uncertainties. Hence, these techniques usually result in conservative solutions. In these techniques, the notion of tubes and invariant sets are used to ensure the safety. For example, in \cite{Tub1} tube-based stable nonlinear controller is designed to keep the trajectory of the system in the tube around the nominal trajectory. To design such controllers, a convex optimization based on contraction stability theory and also lossless convexification are provided. In \cite{Tube2} ellipsoidal approximations of invariant sets along the given nominal trajectory are obtained for nonlinear systems with bounded uncertainty sets. %They provide a bi-level NMPC optimization strategy to optimize the nominal trajectory and approximate the invariant sets. 
In \cite{Russ2}, iterative semidefinite programs based on sums of squares optimization and polynomial Lyapunov functions are provided to construct the tubes and controllers for control-affine polynomial systems.

Probabilistic approaches consider probability distributions of the uncertainties to satisfy the control objectives and probabilistic safety constraints.
However, these approaches are limited to a particular class of dynamical systems and probabilistic uncertainties. For example, in (\cite{Bool1,Bool2}) convex linear program for control of linear systems in the presence of additive Gaussian uncertainties is provided. In this approach, control objectives and safety constraints are written in terms of the mean and covariance of Gaussian distributions of the states of the system. In \cite{Russ4} polynomial systems with additive Wiener processes are addressed. For a given controller, the upper bound of the probability that states of the system leave a finite region of state space over a finite time is obtained. For this purpose, classical super martingale results are leveraged and semidefinite programs to search for exponential barrier functions are provided. 

Sampling-based approaches such as Monte Carlo based techniques (\cite{Samp1,Samp2}) look for controllers that satisfy control objectives for the all sampled uncertainties. Being a randomized approach, no analytical guarantees can be provided. In this paper, we consider probabilistic polynomial dynamical systems and leverage the theory of measures and moments to design polynomial feedback controllers. Measures have been used for safety verification and also control of deterministic polynomial systems. For example, in (\cite{ri1},\cite{ri2}) measures are used for risk estimation in uncertain environments. In \cite{Did1}, notion of occupation measures are used for robustness analysis of given controllers. In \cite{Did2}, occupation measures are used to find the inner approximations of region of attractions sets. In \cite{ Lass2}, optimal control of polynomial systems are addressed and occupation measure based semidefinite programs are provided to find open-loop control inputs. In \cite{Russ6}, occupation measure based feedback controllers are designed to maximize the backward reachable sets of control-affine polynomial systems. This technique, in \cite{Russ7}, is extended to design stabilizing controllers for hybrid polynomial systems. 

In (\cite{Ashk0, Ashk00}), we provide time-invariant polynomial feedback and also model predictive controllers to derive the states of polynomial systems to the goal set in the presence of bounded probabilistic uncertainties. Provided techniques rely on measure and moment based semidefinite program formulation of chance constrained optimization problems (\cite{Ashk,Ashk2,Ashk3}). In the provided controller design approaches, complexity and size of the optimization depend on the length of the planning horizon. Hence, such formulations are limited to the short planning horizons.% (approximately 3 planning time steps).

In this paper, we leverage the measure and moments theory to address control of probabilistic nonlinear systems given nominal trajectories. We formulate this problem as a chance optimization problem where we maximize the probability that states of the system follow the nominal trajectory and remain in the neighborhood of the nominal trajectory despite all uncertainties. 
To address the control problems with long planning horizons, instead of solving a large chance optimization over the planing horizon, we provide a sequence of smaller chance optimization problems. In the provided procedure at each time step, we i) propagate the probability distribution of uncertainties through the nonlinear system to find the probability distributions of the states of the system and ii) solve a small chance optimization problem to obtain the feedback gains. To solve such a optimization problem, we leverage results on chance optimization problems. In (\cite{Ashk,Ashk2,Ashk3}), building on the theory of measure and moments, we provide convex relaxations to efficiently solve chance optimization problems.  Provided approach in this paper deals with bounded and unbounded probabilistic uncertainties.

The outline of the paper is as follows: In Section 2, we cover the preliminary results on polynomials, measures and moments; Section 3 includes the nonlinear controller design problem formulation; Section 4 provides a sequence of chance optimization problems to solve the controller design problem; in Section 5, we provide convex relaxations to solve the chance optimization problems and design the controllers; in Section 6, we address the uncertainty propagation problem to obtain the probability distribution of the states; in Section 7, we present numerical examples, followed by some concluding remarks in Section 8.

%%%%%%%%%%%%%%%%%%%%%%%%%%%%%%%%%%%%%%%%%%%%%%%%%%%%%%%%%%%%%%%%%%%%%%%%%%%%%%%%%%%%%%%%%

\section{ Notation and Preliminary Results}\label{Notation}

%\subsection{ Notations and Definitions}
\label{sec:definitions}

This section covers some basic definitions of polynomials, measures and moments (\cite{SOS1,SOS2,Ashk}). 

\textbf{Polynomials:} Let {\small $\mathbb{R}[x]$} be the set of real polynomials in the variables {\small $x \in \mathbb{R}^n$}. Given {\small $\cP\in\mathbb{R}[x]$}, we represent {\small $\cP$} as {\small $\sum_{\alpha\in\mathbb{N}^n} p_\alpha x^\alpha$} using the standard basis {\small $\{x^\alpha\}_{\alpha\in \mathbb{N}^n}$} of {\small$\mathbb{R}[x]$}, and {\small $\mathbf{p}=\{p_\alpha\}_{\alpha\in\mathbb{N}^n}$} denotes the polynomial coefficients and {\small$\alpha \in \mathbb N ^n$}, e.g., {\small ${\alpha}=(\alpha_1,...,\alpha_n)$} with {\small$\alpha_i$} in {\small$\mathbb{N}$}. Also, {\small$\mathbb R_{\rm d}[x] \subset \mathbb R [x]$} denotes the set of polynomials of degree at most {\small$d\in \mathbb{N}$}. For example a polynomial of degree at most {\small$d=2$} in {\small$x_1$} and {\small$x_2$} ({\small$n=2$}) can be represented as {\small $\cP= \sum_{\alpha={(\alpha_1,\alpha_2)}} p_{\alpha_1\alpha_2} x_1^{\alpha_1}x_2^{\alpha_2}  \in\mathbb R_{\rm 2}[x]$} where {\footnotesize $\alpha={(\alpha_1,\alpha_2)} \in \{(0,0),(1,0),(0,1),(2,0),(1,1),(0,2)\}$} and  {\small$p_{\alpha_1\alpha_2}$} are the coefficients. In this paper, we use polynomials to represent dynamical systems, control inputs, and flow-tubes.

\textbf{Measures:} 
Nonnegative measure {\small$\mu$} is a function that assigns a nonnegative real value to the sets. For example, probability measure (probability distribution) {\small$\mu$} is a function from the set of events to {\small$[0,1]$}, e.g., {\small$\mu(A)=\int_{A} d\mu$}. As another example, Lebesgue measure {\small$\mu_{Leb}$} measures the size of the given set, e.g., {\small$\mu_{Leb}(A)=\int_{A} d\mu_{Leb}= \int_{A} dx$}. Support of the measure {\small$\mu$} denoted by {\small$supp(\mu)$} is the smallest closed set that contains all the sets with nonzero measure. For example support of the uniform probability distribution defined on the interval {\small$[0,1]$} is defined as  {\small$supp(\mu)=[0,1]$}. In this paper, we use measures and probability distributions to obtain equivalent convex chance optimization problem.

\textbf{Moments:} Given a measure {\small$\mu$} supported in {\small$\mathbb{R} ^n$} and {\small$\alpha \in \mathbb N ^n$}, e.g., {\small ${\alpha}=(\alpha_1,...,\alpha_n)$} with {\small$\alpha_i$} in $\mathbb{N}$, the moment of order {\small$\alpha$}  of {\small$\mu$} is defined as {\footnotesize$y_{\alpha_1,\alpha_2,...,\alpha_n}=E[x_1^{\alpha_1}x_2^{\alpha_2}...x_n^{\alpha_n}]=\int x_1^{\alpha_1}x_2^{\alpha_2}...x_n^{\alpha_n} d\mu$}. For example, mean and variance of a probability distribution in {\small$\mathbb{R}$} can be represented in terms of the moments as $E[x]=y_1, E[(x-E(x))^2]=y_2-y_1^2$. Also, for a measure in {\small$\mathbb{R}^2$}, the moment sequence up to order {\small$\alpha=2$} is as {\small $\mathbf{y}=(y_{00},y_{10},y_{01},y_{20},y_{11},y_{02})$}. 
In this paper, we use the moment representation of the probability distributions to obtain relaxed convex chance optimization problem. For this purpose, we need the following matrices and lemmas.

\textbf{Moment Matrix:} Given {\small$d\geq 1$} and a sequence of moments up to order {\small$2d$} denoted by {\small$\mathbf{y}^{2d}$} the moment matrix {\small$M_d(\mathbf{y}^{2d})$} is a symmetric matrix of the form {\small$M_d(\mathbf{y}^{2d})=E[\mathcal{B}'_d\mathcal{B}_d]$} where, {\small$\mathcal{B}_d$} is the vector of monomial basis up to order {\small$d$}. Moment matrix {\small$M_d(\mathbf{y}^{2d})$} contains all the moments up to order $2d$. For instance, let $d=2$ and $n=2$ and {\small$\mathcal{B}_2=[1,x_1,x_2,x_1^2,x_1x_2,x_2^2]$}; then the moment matrix {\small$M_2(\mathbf{y}^{4})=E[\mathcal{B}'_2\mathcal{B}_2]$} reads as 
 
{\tiny \begin{equation} \label{moment matrix exa}
	M_2\left({\mathbf y^4}\right)=\left[ \begin{array}{c}
	
	\begin{array}{ccc} y_{00} \ | & y_{10} & y_{01}| \end{array}
	\begin{array}{ccc} y_{20} & y_{11} & y_{02} \end{array}
	\\
	\begin{array}{ccc} - & - & - \end{array}
	\ \ \ \  \begin{array}{ccc} - & - & - \end{array}
	\\
	
	\begin{array}{ccc} y_{10}\ | & y_{20} & y_{11}| \end{array}
	\ \begin{array}{ccc} y_{30} & y_{21} & y_{12} \end{array}
	\\
	
	\begin{array}{ccc} y_{01}\ | & y_{11} & y_{02}| \end{array}
	\ \begin{array}{ccc} y_{21} & y_{12} & y_{03} \end{array}
	\\
	
	\begin{array}{ccc} - & - & - \end{array}
	\ \ \ \ \  \begin{array}{ccc} - & - & - \end{array}
	\\
	
	\begin{array}{ccc} y_{20}\ | & y_{30} & y_{21}| \end{array}
	\ \begin{array}{ccc} y_{40} & y_{31} & y_{22} \end{array}
	\\
	
	\begin{array}{ccc} y_{11}\ | & y_{21} & y_{12}| \end{array}
	\ \begin{array}{ccc} y_{31} & y_{22} & y_{13} \end{array}
	\\
	
	\begin{array}{ccc}y_{02}\ | & y_{12} & y_{03}| \end{array}
	\ \begin{array}{ccc} y_{22} & y_{13} & y_{04} \end{array}
	
	\end{array}
	\right]
	\end{equation}}

%Given $d\geq 1$ and the sequence of moments up to order $2d$ dented by $\mathbf{y}^{2d}$ the moment matrix
%$M_d(\mathbf{y}^{2d})$ is a symmetric matrix of the form $M_d(\mathbf{y}^{2d})=E[\mathcal{B}'_d\mathcal{B}_d]$ where, $\mathcal{B}_d$ is the vector of monomial basis up to order $d$. Moment matrix $M_d(\mathbf{y}^{2d})$ contains all the moments up to order $2d$. For instance, let $d=2$ and $n=2$ and $\mathcal{B}_2=[1,x_1,x_2,x_1^2,x_1x_2,x_2^2]$; then the moment matrix $M_2(\mathbf{y}^{4})=E[\mathcal{B}'_2\mathcal{B}_2]$ reads as

\textbf{Localizing Matrix:} Given a polynomial {\small$\mathcal{P}$} with order of {\small$\delta$} and a sequence of moments up to order {\small$2d$}, localizing matrix is a symmetric matrix of the form {\small$M_{d-r}(\mathbf{y}^{2d},\mathcal{P})=E[\mathcal{P}\mathcal{B}'_{d-r}\mathcal{B}_{d-r}]$} where {\small$r:=\left\lceil\frac{\delta}{2}\right\rceil$}. For example, given polynomial {\small$\mathcal{P}(x_1,x_2)= bx_1-cx^2_2$}, and moments up to order {\small$4$}, the localizing matrix {\small$M_{2-1}(\mathbf{y}^{4},\mathcal{P})=E[\mathcal{P}\mathcal{B}'_{1}\mathcal{B}_{1}]$} reads as:

{\tiny \begin{equation}
	M_1(\mathbf{y}^4;\mathcal P)= 
	\left[ \begin{array}{ccc}
	by_{10}-cy_{02} & by_{20}-cy_{12} & by_{11}-cy_{03} \\
	by_{20}-cy_{12} & by_{30}-cy_{22} & by_{21}-cy_{13} \\
	by_{11}-cy_{03} & by_{21}-cy_{13} & by_{12}-cy_{04} \end{array}
	\right]
	\end{equation}}
We use the following lemma that gives a necessary and sufficient condition for a sequence $\mathbf y$ to have a representing measure $\mu$ supported on the compact semialgebraic set {\small $\mathcal{K}=\{ x \in \reals^n: \mathcal{P}_i(x) \geq 0, i=1,...,n \}$} (\cite{SOS1,SOS2,Ashk}).
\begin{lemma}
	\label{sec2:lem5}
	The sequence {\small$\mathbf y = \{ y_\alpha\}_{\alpha\in\mathbb{N}^{n}}$} is a moment sequence of a measure $\mu$ supported on the compact set $\mathcal{K}$, if and only if
	{\small 
		$M_d(\mathbf y^{2d})\succcurlyeq 0,\quad M_d(\mathbf y^{2d};\mathcal{P}_i)\succcurlyeq 0, i=1,...,n,  d \in \mathbb N.$
	}
\end{lemma}
In Lemma 1, "{\small$\succcurlyeq 0$}" denotes positive semidefinite matrix. Using Lemma 1, we can reformulate linear programs in measures as semidefinite programs in moments, (\cite{Ashk,Ashk2}).
%%%%%%%%%%%%%%%%%%%%%%%%%%%%%%%%%%%%%%%%%%%%

\section{Problem Statement}\label{Formu}

Consider the following uncertain nonlinear system as
\begin{small}
\begin{equation} \label{sys}
\mathbf{x}(k+1)=f(\mathbf{x}(k),\mathbf{u}(k),\mathbf{\omega}(k))
\end{equation}
\end{small}
where, {\small$\mathbf{x}(k)=[x_{1}(k),...,x_{n}(k)]' \in \chi \subset \reals^n$} are the states and {\small$\mathbf{u}(k)=[u_{1}(k),...,u_{m}(k)]' \in \mathcal{U} \subset \reals^m$} are the control inputs of the dynamical system at time step $k$. Also, {\small$\mathbf{\omega}(k) =[\omega_1(k),...,\omega_{l}(k)] \in \Omega \subset \reals^l$} are uncertain parameters such as uncertain model parameters and disturbances at time step $k$ with known probability distribution {\small$p_{\omega_k}(\omega)$}. Also, $f:\mathbb{R}^{n+m+l}\rightarrow \mathbb{R}^n$ is a polynomial vector function that describes the dynamics of the system.

Moreover, initial states {\small$x(0)$} are uncertain with known probability distribution {\small$p_{\mathbf{x}_0}(\mathbf{x})$}. Let {\small$\mathbf{x}^*:=\{\mathbf{x}^*(k), \ k=1,...,T\}$} be the given nominal trajectory and  {\small$\mathbf{u}^*:=\{\mathbf{u}^*(k), \ k=0,...,T-1\}$} be the given nominal open-loop control input. In the absence of uncertainties, system \eqref{sys} follows the given nominal trajectory {\small$\mathbf{x}^*$} under the nominal open-loop control {\small$\mathbf{u}^*$}. \\

\textbf{Flow-Tube:} we define a flow-tube {\small$\mathcal{FT}(k)$} as a neighborhood around the nominal trajectory {\small$\mathbf{x}^*$}, i.e., {\small$\mathbf{x}^*(k) \in \mathcal{FT}(k), \ k=1,...,T$}. More precisely, we represent such neighborhood at each time step $k$ as the following semi-algebraic set described by level sets of polynomials:

\begin{small}
\begin{equation}\label{FT}
\mathcal{FT}(k)=\{ \mathbf{x} \in \reals^n: \ \mathcal{P}_{k_j}(x) \geq 0, \ j=1,...,\ell \}
\end{equation}
\end{small}
where, {\small$\mathcal{P}_{k_j}(x): \reals^n \rightarrow \reals, \ j=1,...,\ell$} are given polynomials.
For example the flow-tube of the form {\small$\mathcal{FT}(k)=\{ \mathbf{x} \in \reals^n: \  ||\mathbf{x}(k)-\mathbf{x}^*(k) ||_2^2 \leq \epsilon(k)  \}
$} where {\small$\epsilon(k) \in \reals^{+}$} represents a ball around the given nominal trajectory {\small$\mathbf{x}^*(k)$}, at each time $k$.\\

\textbf{Time-Varying Feedback Controller:} Let {\small$\mathbf{\bar{x}}(k)=\mathbf{x}(k)-\mathbf{x}^*(k)$} be the error state vector. We look for the control input of the form {\small$\mathbf{u}(k)=\mathbf{\bar{u}}(k)+\mathbf{u}^*(k)$} where {\small$\mathbf{\bar{u}}(k)$} is polynomial feedback control in {\small$\mathbf{\bar{x}}(k)$}. 
More precisely, {\small$\mathbf{u}(k)=[u_{1}(k),...,u_{m}(k)]' \in \mathcal{U} \subset \reals^m$} 
are polynomial feedback control inputs of the form

\begin{small}
\begin{align}\label{Con1}
& \mathbf{u}(k)=\mathbf{\bar{u}}(k)+\mathbf{u}^*(k) \\
& \bar{u}_{i}(k)= \sum_{\alpha\in\mathbb{N}^n} {g_\alpha}_i(k) \mathbf{\bar{x}}(k)^\alpha \  \in\mathbb{R}_d[x],\ i=1,...,m \subeqn
\end{align}
\end{small}
where, {\small$\{\mathbf{\bar{x}}(k)^\alpha\}_{\alpha\in \mathbb{N}^n}$} are the standard monomials and {\small$\mathbf{G}_{i}(k)=[{g_{\alpha}}_i(k), \alpha\in\mathbb{N}^n]$} are the polynomial coefficients vector that represents the feedback gains of input $i$ at time step $k$. Also, the feedback control in \eqref{Con1} should satisfy the control and gain constraints as follows:

\begin{small}
	\begin{align}\label{Con2}
& \mathbf{u}(k) \in \mathcal{U}=[a_1,b_1]\times...\times[a_m,b_m]  \\
& {g_{\alpha}}_i(k)  \in [{a_i}_{\alpha},{b_i}_{\alpha}], \ i=1,...,m, \alpha \in \mathbb{N}^n \subeqn
\end{align}
\end{small}
where, {\small$a_i, b_i \in \reals$}, for $i=1,...,m$, and {\small${a_i}_{\alpha}, {b_i}_{\alpha} \in \reals$}, for {\small$i=1,...,m$, $\alpha \in \mathbb{N}^n$}.
Now, we define the controller design problem as follows.

\textbf{Controller Design Problem:} Given the probabilistic nonlinear system in \eqref{sys}, nominal trajectory {\small$\{\mathbf{x}^*, \mathbf{u}^*\}$}, and associated flow-tube {\small$\mathcal{FT}$} in \eqref{FT}, we aim to design a nonlinear time-varying state feedback of the form \eqref{Con1} satisfying the constraints \eqref{Con2} such that states of the system follow the nominal trajectory and remain in the given flow-tube, despite all uncertainties, i.e., {\small$\mathbf{x}(k) \in \mathcal{FT}(k), \ k=1,...,T$}. Due to the probabilistic uncertainties, we formulate the controller design problem as a chance optimization problem where we aim to maximize the probability of satisfying control objectives in presence of probabilistic uncertainties.
%that states of the system remains inside the flow-tube in presence of probabilistic uncertainties. 
More precisely, we define the following problem.

\textbf{Chance Optimization Problem: } Find the polynomial controller gains {\small $\mathbf{G}_{i}(k)=[{g_{\alpha}}_i(k), {\alpha\in\mathbb{N}^n}], i=1,...,m, k=0,...,T-1$}, to maximize the probability that trajectory of system {\small $\{\mathbf{x}(k), k=1,...,T \}$} remains inside the given flow-tube {\small $\{\mathcal{FT}(k), k=1,...,T \}$}
considering dynamical and input constraints \eqref{sys}, \eqref{Con1}, \eqref{Con2}
 by solving the following optimization problem:

\begin{footnotesize}
\begin{align}
		&\max_{\mathbf{G}_i(k)|_{i=1,k=0}^{i=m,k=T-1}} \hbox{Probability}_{p_{\mathbf{x}_0},p_{\mathbf{\omega}_k}|_{k=0}^{T-1}}( \cap_{k=1}^{T} \{\mathbf{x}(k) \in \mathcal{FT}(k)\})
		\label{P1}\\
	&	\hbox{s.t.}\quad  \mathbf{G}_{i}(k)=[{g_{\alpha}}_i(k), \alpha\in\mathbb{N}^n], i=1,...,m \subeqn\\
		&  \mathbf{x}(k+1)=f(\mathbf{x}(k),\mathbf{u}(k),\omega(k)) \subeqn\\
		& \mathbf{u}(k)=\mathbf{\bar{u}}(k)+\mathbf{u}^*(k) \subeqn\\
		& \bar{u}_{i}(k)= \sum_{\alpha\in\mathbb{N}^n} {g_\alpha}_i(k) \mathbf{\bar{x}}(k)^\alpha, i=1,...,m \subeqn \\
		& \mathbf{\bar{x}}(k)=\mathbf{x}(k)-\mathbf{x}^*(k) \subeqn\\
		& \mathbf{u}(k) \in \mathcal{U}=[a_1,b_1]\times...\times[a_m,b_m]  \subeqn\\
		&  {g_\alpha}_i(k)  \in [{a_{\alpha}}_i,{b_{\alpha}}_i], \ i=1,...,m, \alpha \in \mathbb{N}^n \subeqn\\
		&  x(0) \sim p_{x_0}(x), \omega(k) \sim p_{\omega_k}(\omega) \subeqn\\
		& k=1,...,T-1 \subeqn
\end{align}
\end{footnotesize}
This problem is nonconvex and computationally hard \cite{Ashk}. Moreover, its complexity increase as the length of planning horizon $T$ increases,\cite{Ashk0}. In the following sections, we will provide a sequence of small chance optimization problems and convex relaxations to efficiently solve the chance optimization in \eqref{P1}.

\textit{Assumption on Uncertainties}.
We assume that probability distribution of uncertainties {\small$p_{\mathbf{x}_0}, p_{\omega_k}, k=0,...,T-1$} are independent.
Also, provided approach works for probabilistic bounded and unbounded uncertainties. In the case of unbounded uncertainties, we assume that probability distributions are completely determined by their moments (\textit{moment determinate}), e.g, normal distribution can completely be determined by its first and second moment, \cite{Lass3}.

\section{Sequential Chance Optimization}
In this section, we provide a sequence of small chance optimization problems to solve the chance optimization \eqref{P1} over the planning horizon {\small $k=1,...,T$}.
Consider the probabilistic nonlinear system in \eqref{sys} and cost function of chance optimization problem in \eqref{P1}. 
According to system \eqref{sys}, states {\small$\mathbf{x}(k+1)$} depend on states {\small$\mathbf{x}(k)$}, uncertain parameter {\small$\mathbf{\omega}(k)$}, and control input {\small$\mathbf{u}(k)$}. Hence, the probability distribution of {\small$\mathbf{x}(k+1)$} denoted by {\small$p_{\mathbf{x}_{k+1}}$} depends only on probability distributions {\small$p_{\mathbf{x}_{k}}$} and {\small$p_{\mathbf{\omega}_{k}}$}, and control input {\small$\mathbf{u}(k)$}. Therefore, we rewrite the cost function {\small $\hbox{Probability}_{p_{\mathbf{x}_0},p_{\mathbf{\omega}_k}|_{k=0}^{T-1}}( \cap_{k=1}^{T} \{\mathbf{x}(k) \in \mathcal{FT}(k)\})$} as {\small $\Pi_{k=0}^{k=T-1} \hbox{Probability}_{p_{\mathbf{x}_{k}},p_{\mathbf{\omega}_{k}}}( \mathbf{x}(k+1) \in \mathcal{FT}(k+1))$} and define the following chance optimization at time {\small$k$}:

\textbf{Sequential Chance Optimization:} At each time {\small $k=0,...,T-1$}, find the polynomial control gains {\small $\mathbf{G}_i(k), i=1,...,m$} to maximize the probability that {\small $x(k+1)$} remains inside the flow-tube {\small $\mathcal{FT}(k+1)$} by solving the following optimization problem: 

 \begin{footnotesize}
	\begin{align}
	\max_{\mathbf{G}_i(k)|_{i=1}^{m}} &\hbox{Probability}_{p_{\mathbf{x}_{k}},p_{\mathbf{\omega}_{k}}}( \mathbf{x}(k+1) \in \mathcal{FT}(k+1))
	\label{P2}\\
		\hbox{s.t.}\quad &  \mathbf{G}_{i}(k)=[{g_{\alpha}}_i(k), \alpha\in\mathbb{N}^n], i=1,...,m \subeqn\\
	&  \mathbf{x}(k+1)=f(\mathbf{x}(k),\mathbf{u}(k),\omega(k)) \subeqn\\
	& \mathbf{u}(k)=\mathbf{\bar{u}}(k)+\mathbf{u}^*(k) \subeqn  \label{P2_ubar0} \\
	& \bar{u}_{i}(k)= \sum_{\alpha\in\mathbb{N}^n} {g_\alpha}_i(k) \mathbf{\bar{x}}(k)^\alpha,\ {\tiny i=1,...,m} \subeqn \label{P2_ubar} \\
	& \mathbf{\bar{x}}(k)=\mathbf{x}(k)-\mathbf{x}^*(k) \subeqn \label{P2_xbar}\\
	& \mathbf{u}(k) \in \mathcal{U}=[a_1,b_1]\times...\times[a_m,b_m]  \subeqn \label{P2_ucon1}\\
	&  {g_\alpha}_i(k)  \in [{a_{\alpha}}_i,{b_{\alpha}}_i], \ i=1,...,m, \alpha \in \mathbb{N}^n \subeqn  \label{P2_ucon2}
	\end{align}
\end{footnotesize}
Note that in chance optimization \eqref{P2}, instead of maximizing the probability in terms of trajectory of the system, we maximize the probability in term of the states at single time step $k$.
Hence, to design a time-varying polynomial controller over {\small$k=0,...,T-1$}, we need to find the probability distribution of states {\small$p_{\mathbf{x}_{k}}$} and solve a chance optimization of the form \eqref{P2} at each time step {\small$k$}. In the next sections, we will provide the convex relaxation of the chance optimization \eqref{P2} and address the uncertainty propagation to obtain the probability distribution of states {\small$p_{\mathbf{x}_{k}}$}.

 \section{Convex Formulation}
In this section, we leverage the recent results on "chance optimization" based on the theory of measure and moments (\cite{Ashk, Ashk2, Ashk3}) to obtain the convex relaxation of the chance optimization in \eqref{P2}. Consider the cost function of chance optimization \eqref{P2}. Using system model \eqref{sys} and flow-tubes \eqref{FT}, we rewrite the cost function as follows

\begin{footnotesize}
\begin{equation}\label{cc1}
\hbox{Probability} \left(  \mathcal{P}_{{k+1}_j} \left( f \left( \mathbf{x}(k),\mathbf{\omega}(k),\mathbf{\bar{u}}(k)+\mathbf{u}^*(k) \right) \right) \geq 0 |_{ j=1}^\ell  \right) 
\end{equation}
\end{footnotesize}
Using \eqref{P2_ubar} and \eqref{P2_xbar}, we represent {\small$\mathbf{\bar{u}}(k)$} in terms of states of the system. %We replace {\small$\mathbf{\bar{u}}(k)$} with \eqref{P2_ubar} and \eqref{P2_xbar} and 
For simplicity, we represent the obtained polynomial $\mathcal{P}_{{k+1}_j}(f(.))$ with polynomial {\small $\hat{\mathcal{P}}_{{k+1}_j}( \mathbf{x}(k),\mathbf{\omega}(k), \mathbf{G}_i(k)|_{i=1}^m)$}. Moreover, since control input {\small$\mathbf{u}(k)$} depends on uncertain states {\small$\mathbf{x}(k)$}, constraint in \eqref{P2_ucon1} is also probabilistic. Hence, we move constraint \eqref{P2_ucon1} inside the probabilistic cost function. Hence, the chance optimization in \eqref{P2} reads as

\begin{footnotesize}
\begin{align}\label{P3}
\max_{\mathbf{G}_i(k)|_{i=1}^{m}} & \hspace{-2mm} \hbox{Probability} \hspace{-1mm} \left(  \hspace{-1mm} \begin{array}{ccc}
	 \hat{\mathcal{P}}_{{k+1}_j}( \mathbf{x}(k),\mathbf{\omega}(k), \mathbf{G}_i(k)|_{i=1}^m) \geq 0 \ |_{j=1}^{\ell}  \\
	   a_i \leq {u}_i(k) \leq b_i \ |_{i=1}^{m} \end{array} \hspace{-1mm} \right)\\
\hbox{s.t.}\quad &  \eqref{P2_ucon2}  \subeqn \label{P3_con}
\end{align}
\end{footnotesize}
In the chance optimization \eqref{P3}, we have two sets of parameters including i) uncertain parameters {\small$\mathbf{x}(k),\mathbf{\omega}(k)$} with probability distributions {\small$p_{\mathbf{x}_{k}},p_{\mathbf{\omega}_{k}}$} and ii) design parameters {\small $\mathbf{G}_i(k)|_{i=1}^m$} that should satisfy the deterministic constraints in \eqref{P3_con}. To obtain a tractable convex relaxation of \eqref{P3}, we take the following steps. \\

\textbf{Equivalent Linear Program:} 
In this step, we obtain a infinite dimensional linear program in terms of probability distributions (\cite{Ashk,Ashk2}). For this purpose, we assign an unknown probability distribution {\small$p_{\mathbf{G}_k}$} to design parameters {\small$\mathbf{G}_i(k), i=1,...,m$}. The support of probability distribution {\small$p_{\mathbf{G}_{k}}$} is in the hyper-cube defined in constraint \eqref{P3_con}, i.e., {\small$\left\lbrace \{ {g_{\alpha}}_i(k)\}_{\alpha\in\mathbb{N}^n}:  {a_{\alpha}}_i \leq {g_\alpha}_i(k) \leq {b_{\alpha}}_i, i=1,...m \right\rbrace$}. Then, we define an equivalent optimization in terms of known probability distributions {\small$p_{\mathbf{x}_{k}},p_{\mathbf{\omega}_{k}}$}
and unknown probability distribution {\small$p_{\mathbf{G}_k}$}. For this purpose, we translate the cost function and constraints of the original chance optimization in \eqref{P3} in terms of probability distributions {\small$p_{\mathbf{x}_{k}},p_{\mathbf{\omega}_{k}}$}, and {\small$p_{\mathbf{G}_k}$}. This results in the following infinite dimensional linear program (Problem 3.2 and Theorem 3.1 in \cite{Ashk},\cite{Ashk2}):  {\small $\max_{p_{\mathbf{G}_k},\mu} \int d\mu, \ \hbox{s.t.} \ \mu \preccurlyeq p_{\mathbf{G}_k} \times p_{\mathbf{x}_{k}} \times p_{\mathbf{\omega}_{k}}$}, {\small $\ \hbox{supp}(\mu) \subset$} {\footnotesize$\left\{ (\mathbf{x}(k),\mathbf{\omega}(k), \mathbf{G}_i(k)|_{i=1}^m):  \begin{array}{ccc}
	 \hat{\mathcal{P}}_{{k+1}_j}( .,.,.) \geq 0 \ |_{j=1}^{\ell}  \\
	   a_i \leq {u}_i(k) \leq b_i \ |_{i=1}^{m} \end{array} \right\}$}, and {\footnotesize$\hbox{supp}(p_{\mathbf{G}_k}) \subset$}{\footnotesize$ \left\lbrace \{ {g_{\alpha}}_i(k)\}_{\alpha\in\mathbb{N}^n}:  {a_{\alpha}}_i \leq {g_\alpha}_i(k) \leq {b_{\alpha}}_i, i=1,...m \right\rbrace$} where $\mu$ is a slack measure. Let, {\small $\{\mathbf{G}^*_i(k), i=1,...,m\}$} and {\small $p^*_{\mathbf{G}_k}$} be the optimal solution of the original chance optimization \eqref{P3} and infinite LP, respectively. Then, the following results holds true: i) the optimal value of infinite LP and original chance optimization are the same, ii) Any {\small $\{\mathbf{G}_i(k), i=1,...,m\} \in supp(p^*_{\mathbf{G}_k})$} is an optimal solution of the original chance optimization, iii) Dirac distribution at {\small $\{\mathbf{G}^*_i(k), i=1,...,m\}$} is an optimal solution of infinite LP, (Problem 3.2 and Theorem 3.1 in \cite{Ashk},\cite{Ashk2}).\\
%%{\footnotesize$\left\{ (\mathbf{x}(k),\mathbf{\omega}(k), \mathbf{G}_i(k)|_{i=1}^m): \hat{\mathcal{P}}_{{k+1}_j}( \mathbf{x}(k),\mathbf{\omega}(k), \mathbf{G}_i(k)|_{i=1}^m) \geq 0|_{j=1}^{\ell} \right\}$},

\textbf{Equivalent Semidefinite Program:} In this step, instead of looking for probability distribution {\small$p_{\mathbf{G}_k}$}, we look for its moment sequences. Let {\small $\mathbf{y}, \mathbf{y}_{\mathbf{G}_k}, \mathbf{y}_{\mathbf{x}_{k}},\mathbf{y}_{\mathbf{\omega}_{k}}$} be the moment sequences of {\small$\mu, p_{\mathbf{G}_k}, p_{\mathbf{x}_{k}}, p_{\mathbf{\omega}_{k}}$}, respectively. Using lemma 2, we translate the cost function and constraints of the linear program in terms of the moment sequence of the distributions. This results in the following infinite dimensional SDP (Problem 3.6 and Lemma 3.2 in \cite{Ashk},\cite{Ashk2}):  {\small $ \sup_{\mathbf{y}, \mathbf{y}_{\mathbf{G}_k}} (\mathbf{y})_\mathbf{0}, \ \hbox{s.t.} $}, {\footnotesize $\ M_{\infty}(\mathbf y)\succcurlyeq 0 ,\ M_{\infty}(\mathbf{y}; \hat{\mathcal{P}}_{{k+1}_j})\succcurlyeq 0, \ j=1,...,\ell,
\ M_{\infty}(\mathbf{y};u_i(k)-a_i)\succcurlyeq 0, \ M_{\infty}(\mathbf{y};b_i-u_i(k))\succcurlyeq 0 \ i=1,...,m$}, {\footnotesize $\left(\mathbf{y}_{\mathbf{G}_k}\right)_\mathbf{0}=1, \ M_{\infty} ({\mathbf y}_{\mathbf{G}_k})\succcurlyeq 0$}, {\footnotesize $M_{\infty}(\mathbf{y}_{\mathbf{G}_k};  {g_\alpha}_i(k)-{a_{\alpha}}_i  )\succcurlyeq 0, M_{\infty}(\mathbf{y}_{\mathbf{G}_k};  {b_{\alpha}}_i-{g_\alpha}_i(k)  )\succcurlyeq 0, \ i=1,...,m, \alpha\in\mathbb{N}^n$}, {\footnotesize $\ M_{\infty} (\mathbf{y}_{\mathbf{G}_k}\times\mathbf{y}_{\mathbf{x}_{k}}\times\mathbf{y}_{\mathbf{\omega}_{k}} -{\mathbf y})\succcurlyeq 0$},
where, {\small $(\mathbf{y})_\mathbf{0}$} is the first element of the moment sequence of measure {\small $\mu$}, and {\small $\mathbf{y}_{\mathbf{G}_k}\times\mathbf{y}_{\mathbf{x}_{k}}\times\mathbf{y}_{\mathbf{\omega}_{k}} -{\mathbf y}$} is the moment sequence of measure {\small $p_{\mathbf{G}_k} \times p_{\mathbf{x}_{k}} \times p_{\mathbf{\omega}_{k}}-\mu$}. The optimal values of the infinite SDP and infinite LP are the same and optimal solution of the infinite SDP is the moment sequence of the optimal solution of the infinite LP.

\textbf{Tractable SDP Relaxation:} In order to obtain a finite SDP, we truncate the infinite dimensional matrices and show that as the size of the matrices increases the optimal solution of the relaxed problem converges to the optimal solution of the original problem (Problem 3.7 and Theorem 3.3 in \cite{Ashk},\cite{Ashk2}). This results in the following SDP relaxation:

{\footnotesize { \begin{align} \label{SDP}
&  \sup_{\mathbf{y}^{2d}, \mathbf{y}^{2d}_{\mathbf{G}_k}} (\mathbf{y}^{2d})_\mathbf{0}, \ \hbox{s.t.} \\
&  M_{d}(\mathbf y^{2d})\succcurlyeq 0,\ M_{d-r_j}(\mathbf{y}^{2d}; \hat{\mathcal{P}}_{{k+1}_j})\succcurlyeq 0, j=1,..,\ell \subeqn \\
& M_{d-r_{u_i}}(\mathbf{y}^{2d};u_i(k)-a_i)\succcurlyeq 0,  M_{d-r_{u_i}}(\mathbf{y}^{2d};b_i-u_i(k))\succcurlyeq 0  |_{i=1}^{m}  \subeqn \\
& \left(\mathbf{y}^{2d}_{\mathbf{G}_k}\right)_\mathbf{0}=1, \ M_{d} ({\mathbf y}^{2d}_{\mathbf{G}_k})\succcurlyeq 0, \subeqn \\
& M_{d-1}(\mathbf{y}^{2d}_{\mathbf{G}_k};  {g_\alpha}_i(k)-{a_{\alpha}}_i  )\succcurlyeq 0, \ i=1,...,m, \alpha\in\mathbb{N}^n, \subeqn \\
& M_{d-1}(\mathbf{y}^{2d}_{\mathbf{G}_k};  {b_{\alpha}}_i-{g_\alpha}_i(k)  )\succcurlyeq 0, \ i=1,...,m, \alpha\in\mathbb{N}^n, \subeqn \\
&  M_{d} (\mathbf{y}^{2d}_{\mathbf{G}_k}\times\mathbf{y}^{2d}_{\mathbf{x}_{k}}\times\mathbf{y}^{2d}_{\mathbf{\omega}_{k}} -{\mathbf y^{2d}})\succcurlyeq 0. \subeqn
	     \end{align}}}
where, {\small $\mathbf{y}^{2d}, \mathbf{y}^{2d}_{\mathbf{G}_k}, \mathbf{y}^{2d}_{\mathbf{x}_{k}},\mathbf{y}^{2d}_{\mathbf{\omega}_{k}}$} are the moment sequences up to order $2d$ of distributions {\small $\mu, p_{\mathbf{G}_k}, p_{\mathbf{x}_{k}}, p_{\mathbf{\omega}_{k}}$}, respectively. {\small $(\mathbf{y}^{2d})_\mathbf{0}$} is the first element of the moment sequence of measure {\small $\mu$}, 
{\small $r_j=\left\lceil\frac{\delta_j}{2}\right\rceil$} where {\small $\delta_j$} is the order of polynomial {\small $\hat{\mathcal{P}}_{{k+1}_j}$}, and {\small$r_{u_i}=\left\lceil\frac{\delta_{u_i}}{2}\right\rceil$} where {\small$\delta_{u_i}$} is the order of polynomial input {\small $u_i(k)$}. Note that, we use moments sequence up to order $2d$ to construct the matrices in SDP \eqref{SDP}.
The following results holds true. i) Optimal value of the finite SDP in \eqref{SDP} is an upper bound of the optimal value of the original chance optimization \eqref{P3} and monotonically converges as $d$ increase, ii) The sequence of optimal solution to the finite SDP converges to the moment sequence of the distributions that are optimal to the infinite LP, (Problem 3.7 and Lemma 3.3 in \cite{Ashk},\cite{Ashk2}).

%To obtain the solution of the original chance optimization in \eqref{P2} from the solution of SDP \eqref{SDP}.
If the solution of the SDP in \eqref{SDP} satisfies the rank condition, we can extract the solution of the original chance optimization in \eqref{P2} solving a linear algebra problem \cite{SOS2}. For example, if the moment matirx has rank one, i.e., $\text{Rank} (M_{d} ({\mathbf y}^{2d}_{\mathbf{G}_k}) )=1$, the optimal solution of the infinite LP in measures {\small $p^*_{\mathbf{G}_k}$}, is a Dirac distribution concentrated on the optimal solution of the original chance optimization \eqref{P3}. Hence, we can approximate the solution of the original chance optimization \eqref{P3} with the first order moments of the optimal solution of the finite SDP, i.e., {\small$\{\mathbf{G}^*_i(k), i=1,...,m\} \approx$} first order moments in the sequence {\small$\mathbf{y}^{*2d}_{\mathbf{G}_k}$}. If the rank condition is not satisfied, we need to increase the SDP relaxation order $d$. For more details and sample code see (\cite{Ashk,Ashk2,Ashk3,SOS2}).

\section{Uncertainty Propagation}
In the previous section, we obtained a convex relaxation of chance optimization \eqref{P2} that relies on the moment information of probability distribution $p_{\mathbf{x}_{k}}$. In this section, we obtain the moment sequence of the probability distribution of states at time $k$. By recursion of the dynamical model in \eqref{sys}, we can write the states $\mathbf{x}(k)$ in terms of the uncertain parameters and control input as follows

\begin{small}
	\begin{equation} \label{m1}
{x_i}(k)={P_{f_k}}_i(\mathbf{x}(0),\mathbf{u}(j)|_{j=0}^{k-1},\mathbf{\omega}(j)|_{j=0}^{k-1}), \ i=1,...,n
\end{equation}
\end{small}
where, ${P_{f_k}}_i(.), i=1,...,n$ are polynomials obtained using dynamical model in \eqref{sys}.
Then, moment of order $\alpha$ reads as
\begin{small}
	\begin{equation} \label{m2}
{y_x}_{\alpha}(k)=E[x_1^{\alpha_1}(k)...x_n^{\alpha_n}(k)]=E[{P_{f_k}}_1^{\alpha_1}(.)...{P_{f_k}}_n^{\alpha_n}(.)]
\end{equation}
\end{small}
Given the control inputs up to time step {\small$k-1$} and moments of uncertainties, {\scriptsize ${y_x}_{\beta}(0)=\hbox{E}[x_1^{\beta_1}(0)...x_n^{\beta_n}(0)]$} and {\scriptsize ${y_{\omega}}_{\gamma}(j)=\hbox{E}[\omega_1^{\gamma_1}(j)...\omega_{l}^{\gamma_l}(j)], j=0,...,k-1$}, and also considering that uncertainties are independent, i.e., {\scriptsize $\hbox{E}[x_1^{\beta_1}(0)...x_n^{\beta_n}(0)\omega_1^{\gamma1}(j)...\omega_{l}^{\gamma_l}(j)\omega_1^{\zeta_1}(j')...\omega_{l}^{\zeta_l}(j')]={y_x}_{\beta}(0){y_{\omega}}_{\gamma}(j){y_{\omega}}_{\zeta}(j')$}, we can rewrite the moments of states at time {\small$k$} in \eqref{m2} as
\begin{small}
 \begin{equation} \label{m3}
{y_x}_{\alpha}(k)= \sum_{j}  c_j\left\{ {y_x}_{\beta_j}(0){y_{\omega}}_{\gamma_j}(0)...{y_{\omega}}_{\eta_j}(k-1) \right\}
\end{equation}
\end{small} 
where, coefficients $c_j$ depends on the parameters of control input and dynamical system. Hence, using \eqref{m3} we can obtain the moments sequence of probability distribution of states at time $k$ in terms of known moments of uncertainties and control inputs up to time step $k-1$.

\section{Implementation and Numerical Results}\label{sec:Res}

\begin{figure}
    \rule[0in]{7.2in}{0.1mm}\\
    {\small \textbf{Algorithm 1:} Control of Probabilistic Nonlinear Systems}\\
    \rule[0.125in]{7.2in}{0.1mm}
    %\vspace{-2mm}

    {\small
    \textbf{Inputs:} probabilistic nonlinear system \eqref{sys}, nominal trajectory $\{\mathbf{x}^*, \mathbf{u}^*\}$, associated flow-tube $\mathcal{FT}$ \eqref{FT}, control input constraints \eqref{Con2}, probability distributions of uncertain parameters {\small$\mathbf{x}(0)$} and {\small$\mathbf{\omega}(k),k=0,...T-1$}.\\
    \textbf{Outputs}:  $\mathbf{G}_{i}(k)=[{g_{\alpha}}_i(k), \alpha\in\mathbb{N}^n]$ feedback gains of inputs over the planning horizon $k=0,...T-1$ in \eqref{Con1}.

    \begin{algorithmic}[1]
    %\STATE $k\gets 0$
    %\STATE $\eta_0\gets0.5~\norm{\grad\gamma\left(x_0\right)-\nu_0 A^{*}\left(\Pi_{C^*}\left(b-A(x_0)\right)\right)}_2$
    \STATE   $k\gets 0$  
    \WHILE{$k\leq T-1$} \label{algeq:stop_ALCC}
    \STATE for non-polynomial dynamical systems, obtain the polynomial approximation using the finite order Taylor expansion of $f$ around the point $\{\mathbf{x^*(k)}, \mathbf{u^*(k)}\}$
    \STATE solve semidefinite program \eqref{SDP} with respect to $p_{\mathbf{x}_k}$
    \STATE construct feedback controller $\mathbf{u}(k)$ in \eqref{Con1}
    \STATE propagate the moments of initial probability distribution $p_{\mathbf{x}_0}$ to the time step $k+1$ using the information of the control inputs $\mathbf{u}(j),\ j=0,...,k$
    \STATE $k\gets k+1$
    \ENDWHILE
    %\RETURN $ x_\ell^{(1)}$
    \end{algorithmic}
    \rule[0in]{7.2in}{0.1mm}
    }
%    \vspace{-0.5in}
%  \caption{Algorithm 1: Control of Probabilistic Nonlinear Systems}\label{tab:1}
\end{figure}
%\vspace{-0.5in}

In this section, two numerical examples are presented that illustrate the performance of the proposed approach. The proposed algorithm to design time varying feedback controller for probabilistic nonlinear systems is described in Algorithm 1. Given a nominal trajectory (maneuver), we use Algorithm 1 to design the controllers in the offline step. Then, we can execute the maneuver by applying the designed controller in real-time. To solve SDP \eqref{SDP}, we use GloptiPoly \cite{Glopti}, which is a MATLAB-based toolbox for moment-based SDP, and Mosek SDP solver. We verify the obtained controller by estimating the probability that trajectory of the system stays in the given flow-tube, i.e., {\small $ \hbox{Probability}_{p_{\mathbf{x}_0},p_{\mathbf{w}_k}|_{k=0}^{T-1}}( \cap_{k=1}^{T} \{\mathbf{x}(k) \in \mathcal{FT}(k)\})$}, using Monte Carlo simulation. For this purpose, we sample from the probability distributions of uncertain parameters {\small$\mathbf{x}(0)$} and {\small$\mathbf{\omega}(k),k=0,...T-1$}, \cite{Ashk}.

%For this purpose, we first uniformly grid $\chi$ into $\bar{N}$ grid-points ($\bar{N}$ depending on the desired precision). Let $\{x^{(i)}\}_{i=1}^{\bar{N}}\subset\chi$ denote the points in the uniform grid. Next, for each grid point $x^{(i)}$, 

%We sample from the probability distributions of uncertain parameters {\small$\mathbf{x}(0)$} and {\small$\mathbf{\omega}(k),k=0,...T-1$}. Let $\{\omega^{(i)}(k)\}_{i=1}^{N_i}$ be $N_i$ i.i.d. sample of random parameter $ \mathbf{\omega}(k),k=0,...T-1$. 
%Then, we approximate {\small $ \hbox{Prob}_{p_{\mathbf{x}_0},p_{\mathbf{w}_k}|_{k=0}^{T-1}}( \cap_{k=1}^{T} \{\mathbf{x}(k) \in \mathcal{FT}(k)\})$}
%by $P^{(i)}_{N_i}:=\frac{1}{N_i}\sum_{k=1}^{N_i}\mathbf{I}\left(x^{(i)},\omega^{(i,k)}\right), \quad $ where 
%$\mathbf{I}\left(x,\omega\right)=1 \ \hbox{if} (x,\omega)\in { (x,\omega): \ \left\lbrace  \mathcal{P}(x,\omega) \leq 0 \right\rbrace }  $, %and $\mathbf{I}\left(x,\omega\right)=0 \ \hbox{otherwise} $. For each $x^{(i)}$, we chose sample size $N_i$ such that $P^{(i)}_{N_i}$ becomes stagnant to further increase in $N_i$. Then, the set of sampled points $x^{(i)}$ whose estimated risk satisfy $ \mathbf{P_{risk}} \geq \epsilon$ represents ${\chi}^{\epsilon}_{risk}$.  \\

\subsection{Example 1: Stabilizing Controller}
In this section, we present a simple example that illustrates the effectiveness of the proposed approach. In this example, we consider the following probabilistic nonlinear system:
{\small$x(k+1)=x(k)+4x(k)^2+0.6x(k)^3+u(k)+0.2\omega(k)-0.1$}, where initial state has normal distribution $x(0) \sim N(0,0.2)$ and uncertainty $\omega$ has triangular distribution over $[0,1]$ with the peak point $0$, i.e., $\omega(k) \sim Tri(0)$. Nominal trajectory and open loop input are given as $x^*(k)=0, k=1,...,8$ and $ u^*(k)=0, k=0,...,7$. The flow-tube is defined as {\small $\mathcal{FT}(k)=\{ x: \  |x| \leq \epsilon(k)\}, \ k=1,...,8$}
where, $\epsilon(1)=0.8, \ \epsilon(2)=0.7, \ \epsilon(3)=0.6, \ \epsilon(4)=0.5, \ \epsilon(5)=0.4, \ \epsilon(6)=0.3, \ \epsilon(7)=0.2, \ \epsilon(8)=0.1$. To ensure that states remains inside the given flow-tube, we design a polynomial feedback of the form {\small $u(k)=g_1(k)x(k)+g_2(k)x(k)^2, \ k=0,...,7$}. Control input and feedback gains should satisfy the following constraints as
{\small $ u(k) \in [-2,2], \ k=0,...,7$} and {\small $ g_1(k),g_2(k) \in [-5,5],\ k=0,...,7$}, respectively.

At time $k=0$, we solve the SDP in \eqref{SDP} considering the moment sequence of probability distributions of $x(0)$ and $\omega(0)$ and obtain control gains $[g_1(0),g_2(0)]$. The $\alpha$-th moment of a Normal distribution with mean $\mu$ and standard deviation $\sigma$ is {\small $y_{\alpha}=\sigma^{\alpha}(-\sqrt{-1}\sqrt{2})^{\alpha}kummerU(\frac{-\alpha}{2}, \frac{1}{2}, \frac{-\mu^2}{2\sigma^2})$} where {\small $kummerU(.,.,.)$} is "confluent hyper-geometric Kummer U function". Also, the $\alpha$-th moment of a triangular distribution $Tri(a)$ over [0,1] and uniform distribution $U[a,b]$ over $[a,b]$ are {\small $y_{\alpha}=\frac{2(1-a^{\alpha+1})}{(\alpha+1)(\alpha+2)(1-a)}$ and  ${y_{\alpha}}=\frac{b^{\alpha+1}-a^{\alpha+1}}{(b-a)(\alpha+1)}$}, respectively.

At time $k=1$, using the obtained control input $u(0)=g_1(0)x(0)+g_2(0)x(0)^2$, we propagate the moments of the probability distribution of $x(0)$ to calculate the moments of the probability distribution of $x(1)$, i.e., {\small ${y_x}_{\alpha}(1)=E[x(1)^{\alpha}]=E[\left( x(0)+4x(0)^2+0.6x(0)^3+u(0)+0.2\omega(0)-0.1\right)^{\alpha}]$}. Then, we solve the SDP in \eqref{SDP} considering the moment sequence of probability distributions of $x(1)$ and $\omega(1)$ and obtain the control gains $[g_1(1),g_2(1)]$. 

%Similarly, at time $k=2$, using the obtained control input $u(1)=g_1(1)x(1)+g_2(1)x(1)^2$, we propagate the moments of the probability distribution of $x(0)$ to obtain the moments of the probability distribution of $x(2)$, i.e., {\footnotesize ${y_x}_{\alpha}(2)=E[x(2)^{\alpha}]$}, where {\footnotesize$x(2)=\left( x(0)+4x(0)^2+0.6x(0)^3+u(0)+0.2\omega(0)-0.1\right) +4\left( x(0)+4x(0)^2+0.6x(0)^3+u(0)+0.2\omega(0)-0.1\right)^2+0.6\left( x(0)+4x(0)^2+0.6x(0)^3+u(0)+0.2\omega(0)-0.1\right)^3+u(1)+0.2\omega(1)-0.1$} and inputs are {\footnotesize $u(0)=g_1(0)x(0)+g_2(0)x(0)^2$} and {\footnotesize $u(1)=g_1(1)\left( x(0)+4x(0)^2+0.6x(0)^3+u(0)+0.2\omega(0)-0.1\right)+g_2(1)\left( x(0)+4x(0)^2+0.6x(0)^3+u(0)+0.2\omega(0)-0.1\right)^2$}. Hence, we can obtain the moments of probability distribution of $x(2)$ completely based on the known moments of probability distributions of $x(0)$, $\omega(0)$, and $\omega(1)$.  Then, we solve the SDP in \eqref{SDP} considering the moment sequence of probability distributions of $x(2)$ and $\omega(2)$ and obtain control gains $[g_1(2),g_2(2)]$. 

\begin{figure}[t!]
	\centering
	\includegraphics[width=8cm, height=5cm]{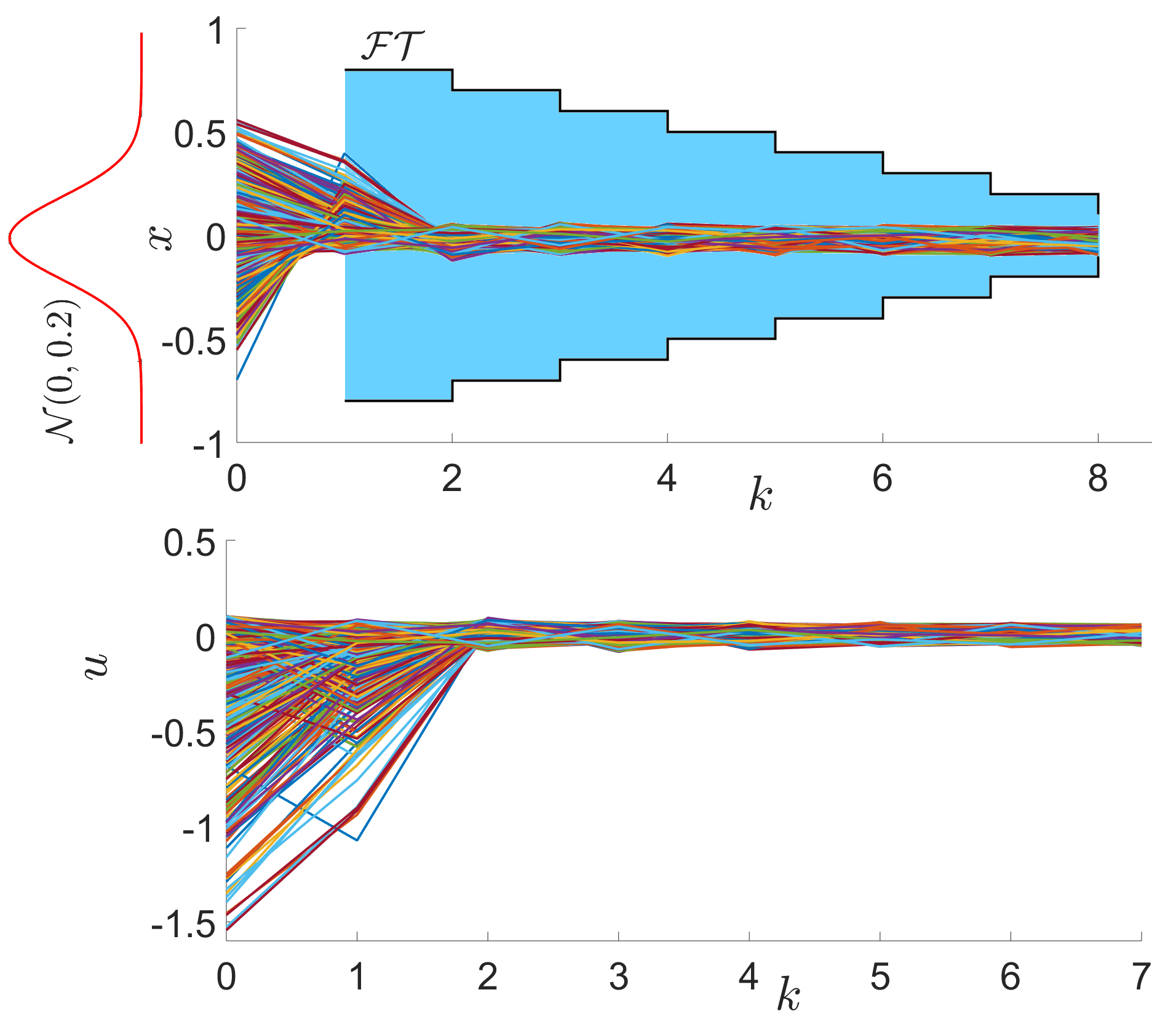}
	\caption{{\footnotesize Example 1- Trajectories of probabilistic nonlinear system under the obtained time-varying polynomial state feedback controller. }}
	\label{fig:1}
\end{figure}

We continue this procedure to obtain all control gains $[g_1(k),g_2(k)], \ k=0,...,7$. The obtained control gains for SDP relaxation order $d=10$ are as {\footnotesize $g_1(k)=[-1.1,   -1.26,   -1.2,   -1.14,   -1.1,   -1.06,   -1.03,   -1] $} and {\footnotesize $ g_2(k)$} {\footnotesize$=[   -2.99,   -3.64,   -3.58,   -3.52 ,  -3.44,   -3.38,   -3.47,  -3.95] $}. Figure \ref{fig:1}, shows the sample trajectories under the obtained controller.% We also, for comparison, design LQR controllers by linearizing the system around the nominal trajectory ( Figure \ref{fig:2}).

%\begin{figure}[t!]
%	\centering
%	\includegraphics[width=8.1cm, height=3cm]{fig/E0/fig22.png}
%	%	\includegraphics[width=8cm, height=3cm]{fig/figObs.png}
%	\caption{{\footnotesize $\mathcal{P}_{risk}(x)$}}
%	\label{fig:1}
%\end{figure}

%
%\begin{figure}[t!]
%	\centering
%	\includegraphics[width=8.1cm, height=3.8cm]{fig/E0/fig3_n.png}
%	%	\includegraphics[width=8cm, height=3cm]{fig/figObs.png}
%	\caption{{\footnotesize Example 1- Results of LQR controller.}}
%	\label{fig:2}
%\end{figure}
%

%\textbf{LQR Results:}\\
%$prob(x(1) \in \mathcal{FT}(1)|x(0)\sim \mathcal{N}(0,0.2))=0.97$ and

%$\{prob(x(k+1) \in \mathcal{FT}(k+1)|x(k) \in \mathcal{FT}(k)), \ k=1,...,7\}=\{0.52,0.55,0.59,0.63,0.69,0.75,0.79\}$

\textbf{Remark 1.} Note that in the presence of unbounded uncertainties, there is always nonzero chance that states of the system leave the flow-tube due to the large disturbances. 
In the provided illustrative example, for $x(0)\geq 0.75$ there exist some $\omega(0)\in [0,1]$ for which trajectory of the system leaves the flow-tube. However, chance of observing such scenarios are close to zero ({\small $Probability(x(0)\geq 0.75)=8.84\times 10^{-5}$}).

\subsection{Example 2: Vehicle Control}
\begin{figure}[t!]
	%\centering
	\includegraphics[width=8.1cm, height=2.5cm]{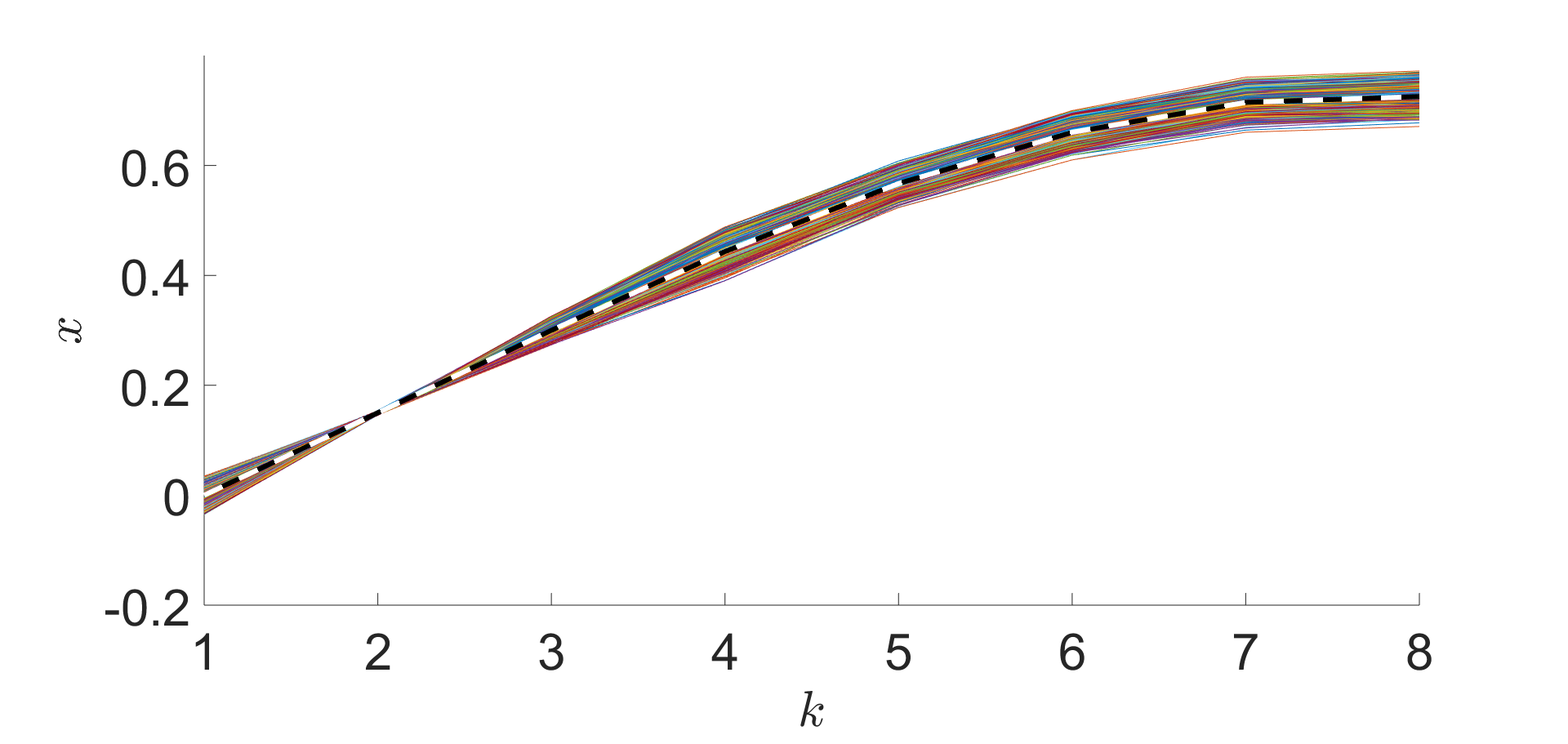}
	\includegraphics[width=8.1cm, height=2.5cm]{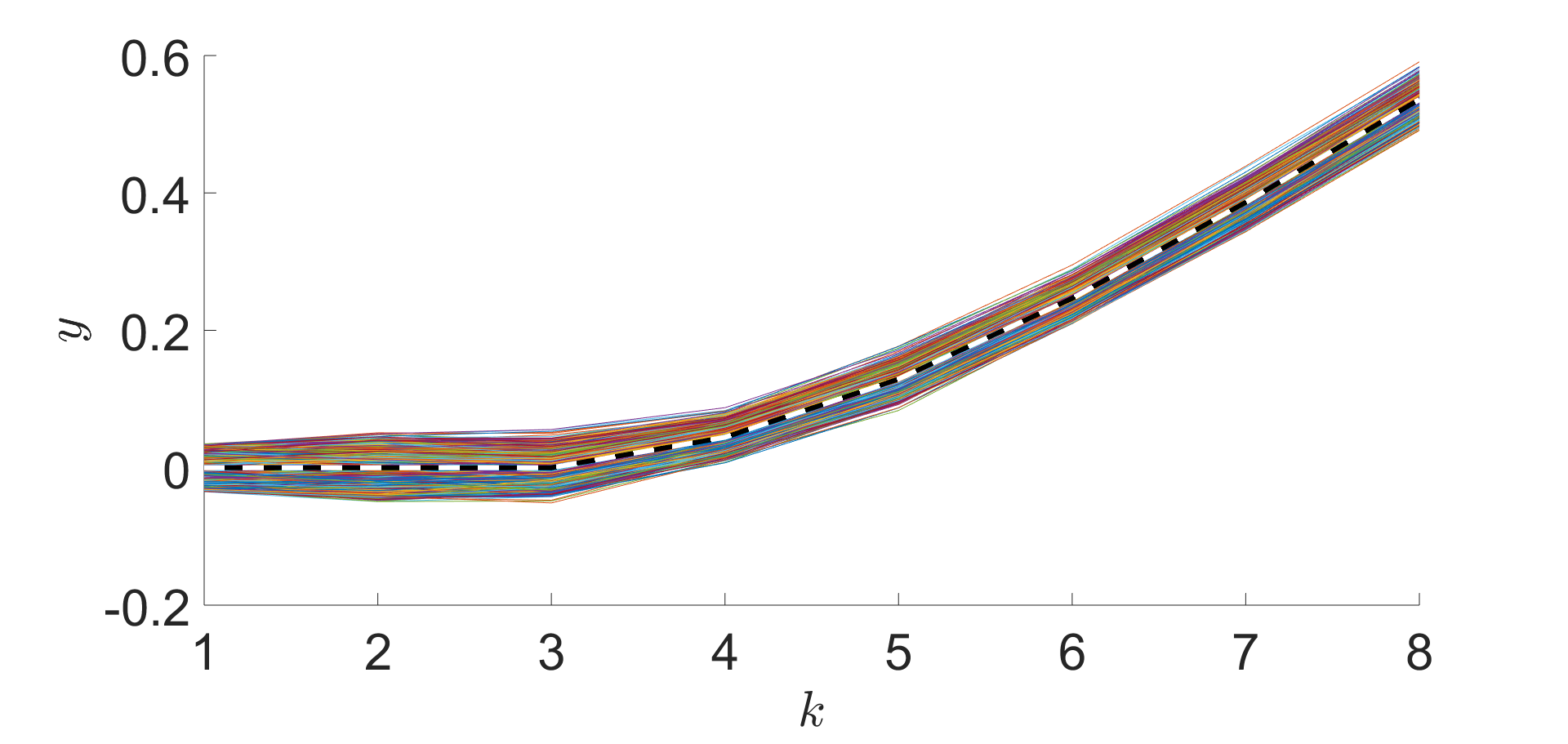}
	\includegraphics[width=8.1cm, height=2.5cm]{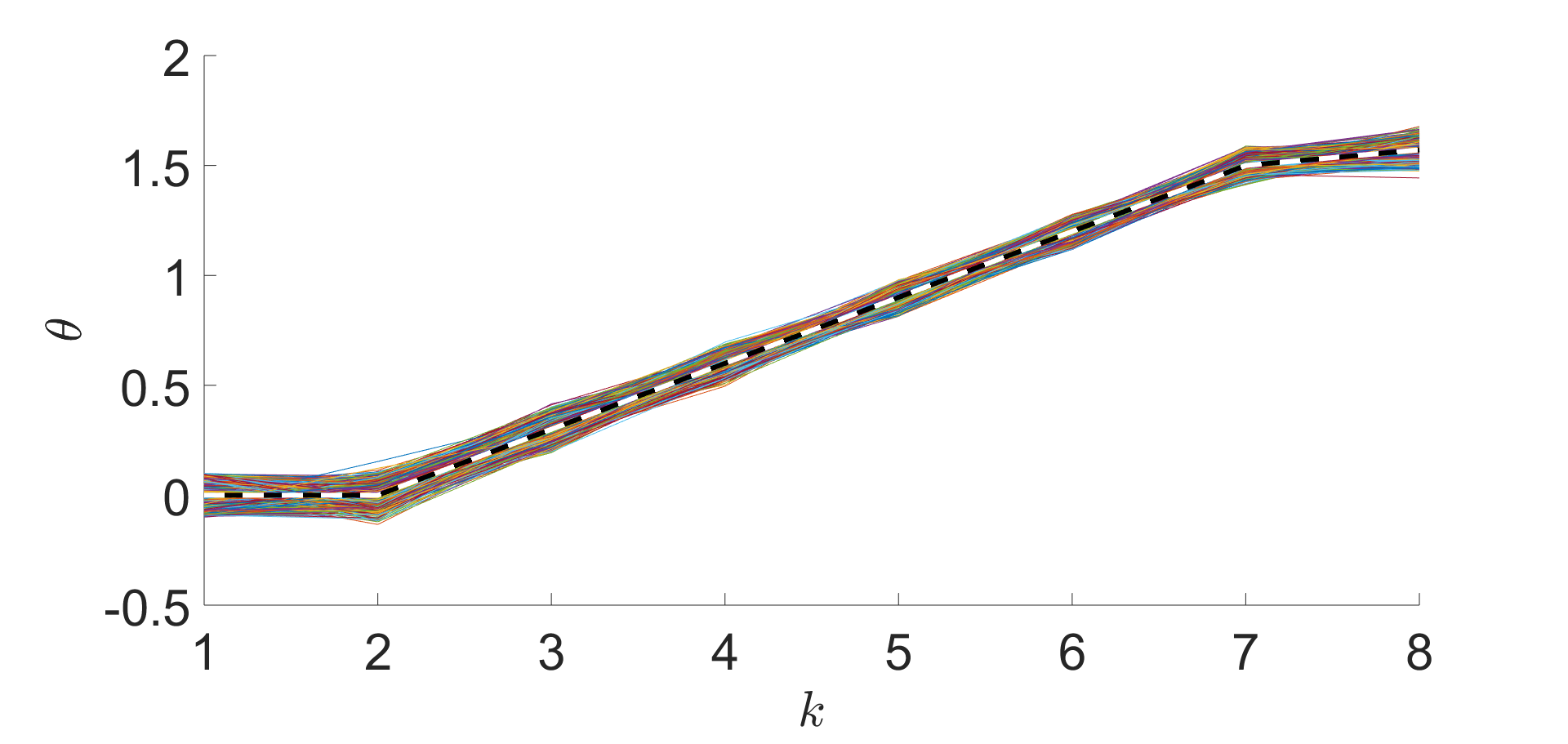}
	\caption{{\footnotesize Sample trajectories of the vehicle under the designed controller}}
	\label{fig:3}
\end{figure}

\begin{figure}[t!]
	\centering
	\includegraphics[width=8cm, height=2.5cm]{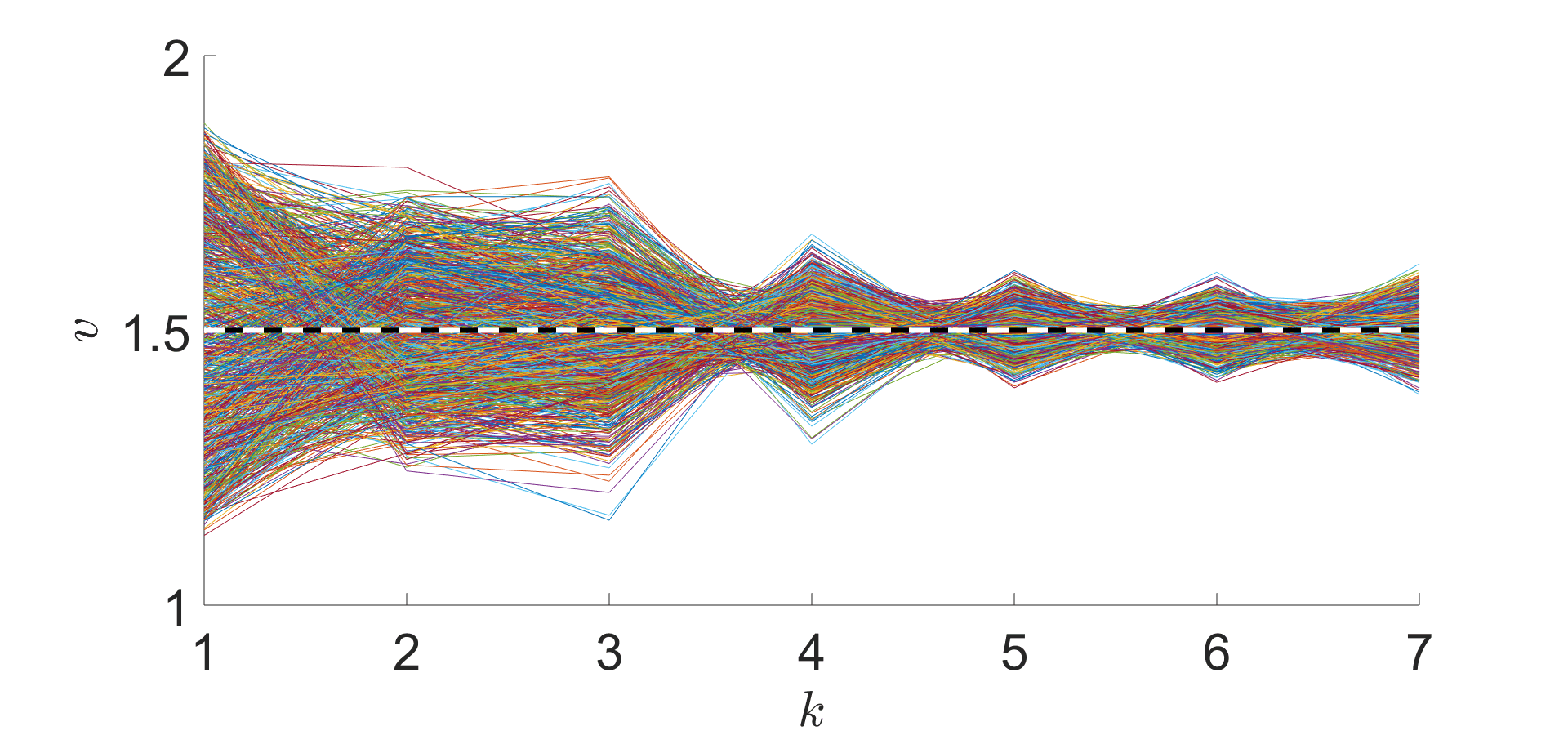}
	\includegraphics[width=8cm, height=2.5cm]{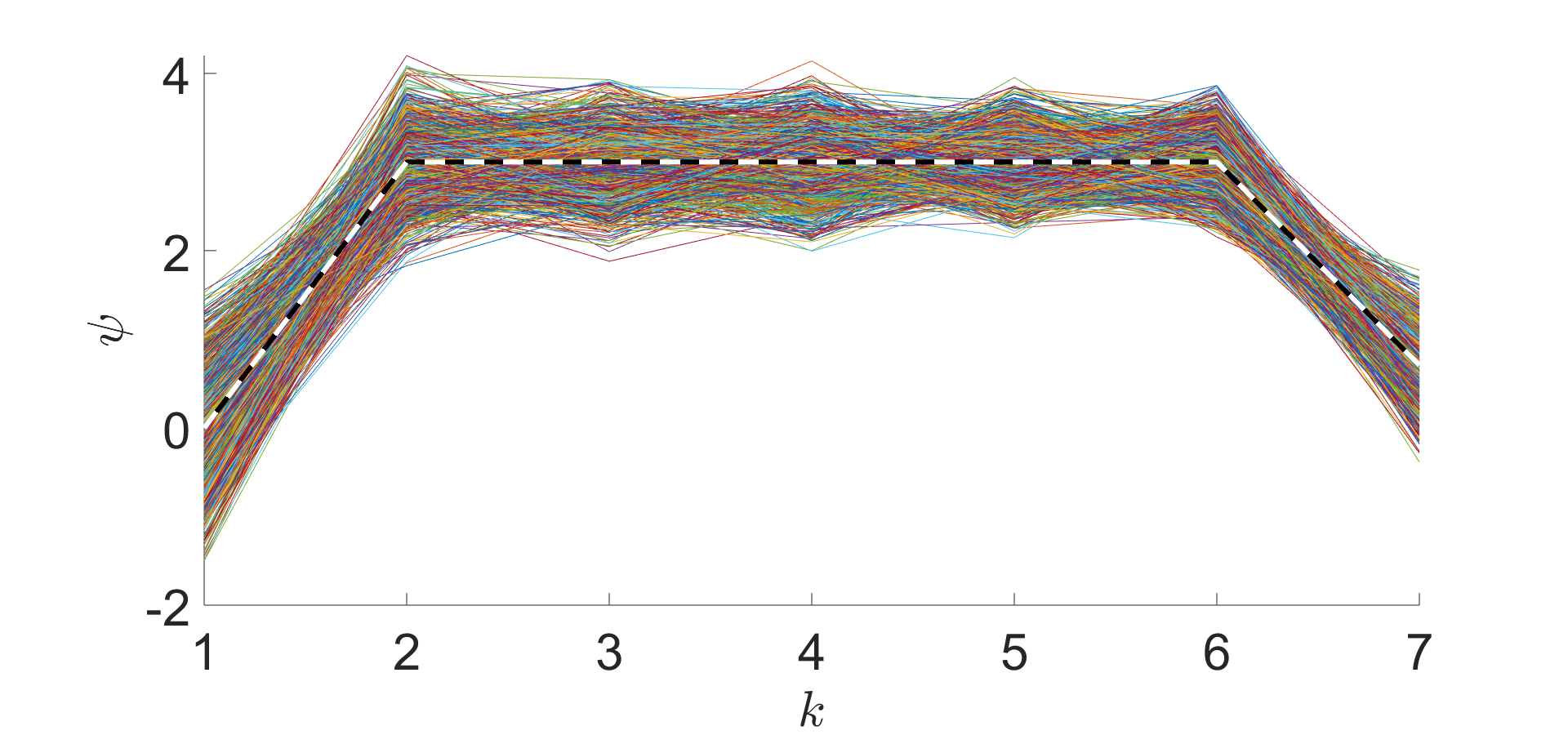}
	\caption{{\footnotesize Sample control inputs of the vehicle}}
	\label{fig:4}
\end{figure}

Dynamics of a race car on a sharp curve is molded as
\begin{footnotesize}
\begin{align}
& x(k+1)=x(k)+\Delta T (v(k)+ 0.1\tilde{v}(k)-0.05)cos(\theta(k))  \notag \\
& y(k+1)=y(k)+\Delta T (v(k)+ \tilde{v}(k))sin(\theta(k))  \notag \\
& \theta(k+1)=\theta(k)+\Delta T ( \psi(k)+0.2\tilde{\psi}(k)-0.1  ) \notag
\end{align}
\end{footnotesize}
where $x,y$ are position, $\theta$ is the steering angle of the vehicle. Inputs are linear velocity $v$ and angular velocity $\psi$. 
Nominal state trajectories and control inputs are given as: {\footnotesize
	 $x^*=\{ 0, 0.15, 0.3, 0.44,0.56,0.66,0.71,0.72 \} $}, {\footnotesize $y^*=\{ 0, 0, 0, 0.04, 0.12, 0.24, 0.38, 0.53 \} $}, {\footnotesize $\theta^*=\{0,0,0.3,0.6, 0.9,1.2,1.5,1.5708\} $} and {\footnotesize $(v,\psi)^*=\{ (1.5,    0),(1.5,    3),(1.5,    3),(1.5,    3),(1.5,    3),(1.5,    3),(1.5,    0.7)\} $}.
The uncertainties such as tire slip are modeled with $0.1\tilde{v}-0.05$ and $0.2\tilde{\psi}-0.1$. At each time step $k$, $\tilde{v}$ and $\tilde{\psi}$ have Beta distribution over $[0,1]$ with parameters $\alpha=4$ and $\beta=4$, i.e., $p_{\tilde{v}_k} = Beta(4,4)$, $p_{{\tilde{\psi}}_k} = Beta(4,4)$. Moment sequence of such Beta distribution are described as {\small $y_{i}=\frac{4+i-1}{8+i-1}y_{i-1}, y_0=1$}.
Moreover, initial states have uniform distributions as {\footnotesize$p_{x_0}=U[-0.07,0.07]$}, {\footnotesize$p_{y_0}=U[-0.07,0.07]$, $p_{{\theta}_0}=U[-0.1,0.01]$}. The flow-tube at each time $k$ is defined as {\footnotesize $\mathcal{FT}(k)=\{ (x,y): \  x^*_k-0.06 \leq x \leq x^*_k+0.06 , y^*_k-0.06 \leq y \leq y^*_k+0.06 \}$}. We are looking for closed loop controllers {\footnotesize $\psi(k)={g_1}_1(\theta(k)-\theta^*(k))+{g_2}_1(x(k)-x^*(k))+{g_3}_1(y(k)-y^*(k))+\psi^*(k)$} and {\small $v(k)={g_1}_2(x(k)-x^*(k))+{g_2}_2(y(k)-y^*(k))+v^*(k)$} to maximize the probability that trajectories of the vehicle remain inside the given flow-tube in the presence of uncertainties. Constraints of the feedback gains and controllers are given as {\footnotesize $-10 \leq {g_1}_1,{g_2}_1,{g_3}_1,{g_1}_2,{g_2}_2 \leq 10$}, {\footnotesize $0 \leq v(k) \leq 2$}, respectively. To obtain the polynomial dynamics, we compute a degree 3 Taylor expansion of the dynamics of the system around the nominal trajectory, at each time $k$. By solving the chance optimizations for relaxation order $d=10$, we obtain the controllers for which vehicle remains in the given flow-tube with probability of one. Obtained results are shown in Figures \eqref{fig:3} and \eqref{fig:4}.

%\begin{figure}[t!]
%	\centering
%	\includegraphics[width=7 cm, height=3.1cm]{fig/E1/Path2.png}
%	\caption{{\footnotesize Flow-tubes and sample trajectories of the vehicle}}
%	\label{fig:5}
%\end{figure}

%%%%%%%%%%%%%%%%%%%%%%%%%%%%%%%%%%%%%%%%%%%%%%%%%%%%%%%%%%%%%%%%%%%%%%%%%%%%%%%%%%%%

%%%%%%%%%%%%%%%%%%%%%%%%%%%%%%%%%%%%%%%%%%%%%%%%%%%%%%%%%%%%%%%%%%%%%%%%%%%%%%%%%%%%%%%%%%%%%

\section{Conclusion}\label{sec:Con}

In this paper, we address the problem of time-varying polynomial state feedback controller design for nonlinear polynomial systems in the presence of bounded and unbounded probabilistic uncertainties. For this purpose, we formulate the controller design problem as a sequence of chance optimization problems where we maximize the probability that states of the uncertain system follow and remain in the given tube around the nominal trajectories. Then, building on the theory of measure and moments, we provide the convex relaxations in the form of semidefinite programs to efficiently solve the obtained chance optimizations. Provided approach deals with both bounded and unbounded probabilistic uncertainties and also long planning horizons. We provide numerical examples on stabilizing controller design and motion planning of uncertain nonlinear systems to illustrate the performance of the proposed approach.

For the future work, we will use the proposed method in control primitive based motion planning where we design library of controllers for different maneuvers in the offline step and execute the right maneuver based on the observed scenarios in real-time.
We will also improve the presented approach to design the controller and, at the same time, smallest flow-tube around the given trajectory to minimize the influence of the uncertainties and deviations from the nominal trajectory. Moreover, we will introduce risk bounded flow-tubes for unbounded uncertainties where the chance of leaving the flow-tube should be bounded by the predefined risk bound.

{}

\end{document}